\documentclass[11pt]{article}
\usepackage{amstext,amsfonts,graphics,latexsym}

\newtheorem{thm}{Theorem} \newtheorem{lemma}[thm]{Lemma}
\newtheorem{lem}[thm]{Lemma}

 \newtheorem{cor}[thm]{Corollary}

\newtheorem{claim}[thm]{Claim}

 \newenvironment{pfof}[1]{\par\vskip2\parsep\noindent{\sc Proof of\ #1. }}{{\hfill $\Box$}
\par\vskip2\parsep}

\newcommand{\ko}{K}

\newcommand{\E}{2}

\newcommand{\ann}{A}
\newcommand{\ret}{R}
\newcommand{\esc}{E}
\newcommand{\Sexc}{S^{\text{exc}}}

\newcommand{\zt}{{\mathbb Z}^2}

\newcommand{\prob}{\mbox{\bf P}}
\newcommand{\expe}{\mbox{\bf E}}
\newcommand{\leb}{{\mathcal L}}

\newcommand{\olle}{H\"{a}ggstr\"{o}m}

\newcommand{\z}{\mathbb Z}
\newcommand{\Z}{\mathbb Z}
\newcommand{\Ex}{\mathbb E}

\newcommand{\N}{\mathbb N}
\newcommand{\dd}{\mathbb D}

\newcommand{\be}{\begin{equation}}
\newcommand{\ee}{\end{equation}}
\newcommand{\0}{{\bf 0}}

\def\qed{\relax\ifmmode\hskip2em \Box\else\unskip\nobreak\hfill
$\Box$\fi}

\newcounter{mycount}
\newenvironment{proof}{\noindent{\sc Proof. }}{{\hfill $\Box$}\par\vskip2\parsep}

\title{A special set of exceptional times for dynamical random walk on $\Z^2$ }
\author{Gideon Amir\thanks{Weizmann Institute, Rehovot, 76100,
Israel. Email: gideon.amir@weizmann.ac.il} \  and Christopher
Hoffman\thanks{Funded in part by an MSRI, DMS grant \#0501102 and
    the University of Washington Royalty Research Fund}}

\begin{document}
\maketitle

\begin{abstract}
Benjamini, \olle, Peres  and Steif introduced the model of
dynamical random walk on $\Z^d$ \cite{ds}.  This is a continuum of
random walks indexed by a parameter $t$.  They proved that for
$d=3,4$ there almost surely exist $t$ such that the random walk at
time $t$ visits the origin infinitely often, but for $d \geq 5$
there almost surely do not exist such $t$.

Hoffman showed that for $d=2$ there almost surely exists $t$ such
that the random walk at time $t$ visits the origin only finitely
many times \cite{H1}.  We refine the results of \cite{H1} for
dynamical random walk on $\z^2$, showing that with probability one
there are times when the origin is visited only a finite number of
times while other points are visited infinitely often.
\end{abstract}


\section{Introduction}

We consider a dynamical simple random walk on $\zt$.  Associated
with each $n$ is a Poisson clock.  When the clock rings the $n$th
move of the random walk is replaced by an independent random
variable with the same distribution. Thus for any fixed $t$ the
distribution of the walks at time $t$ is that of simple random
walk on $\zt$ and is almost surely recurrent.

More formally let $\{Y_n^m\}_{m,n \in \N}$ be uniformly distributed
i.i.d.\ random variables chosen from the set
$\{(0,1),(0,-1),(1,0),(-1,0)\}$. \\Let $\{\tau_n^{(m)}\}_{m \geq 0,n
\in \N}$ be an independent Poisson process of rate one and \\
$\tau_n^{(0)}=0$ for each $n$. Define
$$X_n(t)=Y_n^m$$
for all $t \in [\tau_n^{(m)},\tau_n^{(m+1)}).$ Let
$$S_n(t)=\sum_{i=1}^{n}X_i(t).$$
Thus for each $t$ the random variables $\{X_n(t)\}_{n\in \N}$ are
i.i.d.

The concept of dynamical random walk (and related dynamical
models) was introduced by Benjamini, \olle, Peres and Steif in
\cite{ds}, where they showed various properties of the simple
random walk to hold for all times almost surely (such properties
are called ``dynamically stable"), while for other properties
there a.s.\ exist exceptional times for which they do not hold
(and are then called ``dynamically sensitive" properties). In
particular, they showed that transience of simple random walk is
dynamically sensitive in $3$ and $4$ dimensions, while in
dimension $5$ or higher it is dynamically stable. For more work on
dynamical models see also \cite{levinone} and \cite{dp}.

In \cite{H1} Hoffman showed that recurrence of simple random walk
in two dimensions is dynamically sensitive, and that the Hausdorff
dimension of the set of exceptional times is a.s. $1$. Furthermore
it was showed that there exist exceptional times where the walk
drifts to infinity with rate close to $\sqrt{n}$.

The following problem was raised by Jeff Steif and Noam Berger
during their stay at MSRI: When the dynamical random walk on
$\Z^2$ is in an exceptional time, returning to the origin only a
finite number of times, is it always because it ``escapes to
infinity'', or are there times where the disc of radius $1$ is hit
infinitely often, but the origin is hit only a finite number of
times.

In this paper we prove that such exceptional times exist, and
furthermore that the Hausdorff dimension of the sets of these
exceptional times is $1$. The techniques used are a refinement of
the techniques used in \cite{H1} to show the existence of
exceptional times where the walk is transient.


The proof is divided into two main parts. First, we define a
sequence of events $SE_k(t)$, dependent on $\{S_n(t)\}_{n\in\N}$
such that $\bigcap SE_k(t)$ implies we are in an exceptional time
of the desired form, and show that the correlation between these
events in two different times $t_1,t_2$ is sufficiently small
(Theorem \ref{fmt}). Then we use this small correlation event, in
a manner quite similar to \cite{H1} to show the existence of a
dimension $1$ set of exceptional times.

\section{Definitions, Notation and Main Results} \label{sec:defs}

We will follow the notation of $\cite{H1}$ as much as possible.
Let $\dd$ denote the discrete circle of radius $1$:
$$\dd= \Big \{(1,0),(-1,0),(0,1),(0,-1)\Big\}.$$

 The set of special exceptional times of the process, $\Sexc$, is the set of all times for which the walk
visits $\dd$ infinitely often but visits $(0,0)$ only fintely
often.  More explicitly we define
 $$\Sexc=\bigg\{t:\  \Big|  \{n:S_n(t) \in \dd\}\Big|=\infty \  \text{ and }  \
        \Big|\{n:S_n(t)=\0\} \Big|< \infty \bigg\}.$$
The main result of this paper is:
\begin{thm}\label{expdim1}
$\prob (\Sexc \neq \emptyset)=1$. Furthermore the Hausdorff
dimension of $\Sexc$ is $1$ a.s.
\end{thm}

We will need the following notation.
 Define a series of stopping times \[s_k = k^{10} 2^{2k^2}\] We
 will divide the walk into segments between two consecutive $s_i$:
 $$M_j = \{i\in \N : s_{j-1}\leq i < s_j\},$$

And define  a super segment to be a union of segments
$$W_k = \bigcup_{2^k \leq j < 2^{k+1}}M_j = \{i \in \N: \ s_{2^k-1} \leq
i < s_{2^{k+1}-1}\}.$$


Define a series of annuli

$$A_k = \left\{x\in \Z^2 \ : \ 2^{k^2} \leq |x| \leq k^{10} 2^{2k^2}\right\}.$$

where $|x|$ denotes the standard  Euclidean norm
$|x|=\sqrt{x_1^2+x_2^2}$.

And define event $G_k$ by
$$G_k(t)=\{S_{s_k}(t) \in A_k \}.$$

Choosing these annuli to be of the right magnitude will prove to
be one of the most useful tools at our disposal.
 We will show that a random
walk will typically satisfy $G_k$ (see Lemmas \ref{gbound} and
\ref{gbigbound}). On the other hand, the location of the walk
inside the annuli $A_j$ in the beginning of the segment $M_j$ will
have only a small influence of the probabilities of the events we
will investigate (see Lemmas \ref{mainlem} and \ref{bothreturn}).
Together this provides ``almost independence"
of events in different segments.\\

We will define three events concerned with hitting the disc on the
dynamical random walk.
 We define the event $R_j(t)$, to
be the event that the walk (at time $t$) hits $\dd$ at some step in
in the segment $M_j$, i.e.
$$R_j(t) := \{\exists m\in M_j \text{ such that } S_m(t)\in \dd\}.$$
 (Note that in \cite{H1}, this denoted the event of hitting the
origin, and not $\dd$, but this difference does not alter any of
the calculations.)

We denote by  $SR_j(t)$ the event that the walk (at time $t$) hits
$\dd$ at some step in the segment $M_j$, but does not hit the origin
at any step of $M_j$:
$$SR_j(t) := \left\{\Big(\exists n\in M_j \text{ such that } S_n(t)\in \dd \Big)
    \cap \Big(\forall m\in M_j \ S_m(t)\neq (0,0)\Big)\right \}.$$
It is easy to see that $SR_j(t) \subset R_j(t)$.

We define the third event, $SE_k(t)$, by
$$SE_k(t) = \left\{ \bigcup_{j=2^k}^{2^{k+1}-1}
     \left( SR_j(t) \cap \left( \bigcap_{i=2^k, i \neq j}^{2^{k+1}}(R_i(t))^c\right)\right)\right\}  \bigcap
   \left\{ \bigcap_{i=2^k}^{2^{k+1}}G_i(t)\right\}.$$
This is the event that the walk never hits the origin in the super
segment $W_k$, and there is exactly one segment $M_j, \ 2^{k} \leq j
< 2^{k},$ where the walk hits $\dd$. Additionally this event
requires the walk ends each segment $M_i$ in the annulus $A_i$.

It is easy to see that if the event $\bigcap_{k\geq M} SE_k(t)$
occurs for some integer $M$ then $t\in \Sexc$.\\ Given $2^k\leq j
< 2^{k+1}$ we write
 $$SE_k^{j}(0) = SE_k(0)\cap SR_{j}(0).$$
For any random walk event $E$ we write $E(0,t)$ for $E(0)\cap
E(t)$. Thus for example, we have
 $$SR_j(0,t) = SR_j(0)\cap SR_j(t).$$

For $j_1,j_2 \in [2^k,2^{k+1})$ we write
$$SE^{j_1,j_2}_k(0,t) = SE_k(0,t)\cap SR_{j_1}(0) \cap SR_{j_2}(t).$$

We will use the following notation for conditional probabilities.
Let
$$\prob^{x,k-1}(*)=\prob\left(* |\
      S_{s_{k-1}}(0)=x \right) $$ and
$$\prob^{x,y,k-1}(*)=\prob\left(* |\
      S_{s_{k-1}}(0)=x \mbox{ and } S_{s_{k-1}}(t)=y\right). $$

To prove Theorem \ref{expdim1} we will first prove the following
theorem.
\begin{thm}\label{fmt}
There is a function $f(t)$ such that for any $M$,
$$\frac{\prob(\bigcap_{1\leq k \leq M}SE_k(0,t))}
            {(\prob\bigcap_{1\leq k \leq M}(SE_k(0)))^2} \leq f(t)$$
and $\int_0^1 f(t) < \infty$.
\end{thm}
The first statement in Theorem \ref{expdim1} follows from Theorem
\ref{fmt} and the second moment method. The second statement in
Theorem \ref{expdim1} will follow from bounds that we will obtain
on the growth of $f(t)$ near zero and by Lemma 5.1 of
\cite{peres}.

Since we will be showing that some events have low correlation
between times $0$ and $t$, some of the bounds we use to prove
Theorem \ref{fmt} will only hold when $k$ is sufficiently large
compared to $\frac{1}{t}$. Therefore, for $t\in (0,1)$, we denote
by $K(t)$ the smallest integer greater than $|\log t| + 1$, and by
$K'(t)$ the smallest integer greater than $\log (|\log t| +1) +
1$. (Here and everywhere we use $\log(n) = \log_2(n)$).

 We end this section by stating that we use $C$ as a generic constant
whose value may increase from line to line and lemma to lemma. In
many of the proofs in the next section we use bounds that only
hold for sufficiently large $k$ (or $j$). This causes no problem
since it will be clear that we can always choose $C$ such that the
lemma is true for all $k$.

\section{Proof of Theorem \ref{fmt}}
This section is divided into four parts. In the first, we
introduce some bounds on the behavior  of two dimensional simple
random walk. These are quite standard, and are brought here
without proof. In the second part we quote, for the sake of
completeness, lemmas from \cite{H1} which we will use later on. In
the third part, which constitutes the bulk of this paper, we prove
estimates on $\prob(SE_k(0)\ | \ SE_{k-1}(0))$ and
$\prob(SE_k(0,t) \ | \ SE_{k-1}(0,t))$.  In the fourth part we use
the these bounds to prove Theorem \ref{fmt}.

\subsection{Two Dimensional Simple Random Walk Lemmas}

The main tool that we use are bounds on the probability that
simple random walk started at $x$ returns to the origin before
exiting the ball of radius $n$ with center at the origin. For
general $x$, this probability is calculated in Proposition 1.6.7
on page 40 of \cite{lawler} in a stronger form than given here.
The estimate for $|x|=1$ comes directly from estimates on the
$2-d$ harmonic potential (Spitzer 76, P12.3 page 124 or
\cite{KS}).\\

Let $\eta$ be the  smallest $m>0$ such that $S_m(0)=\0$ or
$|S_m(0)| \geq n$.
\begin{lemma} \label{lawler}
\begin{enumerate}
\item
 There exists $C>0$ such that for all $x$ with $0<|x|<n$
$$\frac{\log(n)-\log|x|-C}{\log(n)}\leq \prob(S_\eta(0)=\0|S_0(0)=x)\leq
\frac{\log(n)-\log|x|+C}{\log(n)}.$$

\item For $|x|=1$ we have the following stronger bound:
$$  P(|S_\eta(0)| \geq n \ | \ S_0(0)=x )= \frac{\pi}{2 \log n + C} + O(n^{-1}). $$
\end{enumerate}
\end{lemma}

We will also use the following standard bounds.
\begin{lemma} \label{leave}
There exists $C>0$ such that for all $x \in \zt$, $n \in \N$ and
$m<\sqrt{n}$
\begin{equation} \label{le}
\prob(\exists n'<n:S_{n'}(0)>m\sqrt{n})\leq \frac{C}{m^2}
\end{equation} and
\begin{equation}\prob\left(|S_{n}(0)-x|<\frac{\sqrt{n}}{m}\right)\leq
\frac{C}{m^2}.\end{equation}
\end{lemma}

Lemma \ref{leave} is most frequently applied in the form of the
following corollary:
\begin{cor}\label{hit0}
There exists $C>0$, such that for any $n_0,N,M>0$, if
$|S_{n_0}(0)|\geq M$ then
$$\prob(\exists n_0 < n < n_0+N \ : \ |S_n(0)| \leq 1) \leq \frac{C N}{M^2}.$$
\end{cor}

\subsection{Some Useful Lemmas From \cite{H1}}

First we list some lemmas from \cite{H1} that we will use in the
proof. ({In \cite{H1} $R_j(0)$ is defined as
$$R_j(0) = \{\exists m\in M_j \text{ such that }
S_m(0)=\0\}.$$ It is easy to see all lemmas in \cite{H1} hold with
our slightly expanded definition.)}

$ $\\
Lemma \ref{gbound} shows that if the walk is in the annulus
$A_{j-1}$ at step $s_{j-1}$, then with high probability the walk
will be in the annulus $A_j$ at step $s_j$. This will allow us to
condition our estimates of various events on the event that
$S_{s_j} \in A_j$.
\begin{lemma} \label{gbound} For any $j$ and $x \in \ann_{j-1}$
$$\prob^{x,j-1} ((G_j(0))^C) \leq \frac{C}{j^{10}}.$$
\end{lemma}

We will usually use the following corollary of Lemma \ref{gbound},
which is just an application of Lemma \ref{gbound} to all segments
inside a super-segment, and follows directly by applying the union
bound:
\begin{cor}\label{gbigbound}
For any $k>0$ and $x \in \ann_{2^k-1}$
\[\prob^{x,2^k-1}\left(\bigg(\bigcap_{2^k \leq l < 2^{k+1}}
G_l(0)\bigg)^c\right)\leq \frac{C}{2^{9k}}.\]
\end{cor}

Lemma \ref{mainlem} shows that we have a good estimate of the
probability of hitting $\dd$ during the segment $M_k$, given the
walk starts the segment inside the annulus $\ann_{k-1}$.
\begin{lemma} \label{mainlem}
There exists $C$ such that for any $k$ and any $x \in \ann_{k-1}$
$$\frac{\E}{k} -\frac {C\log k}{k^2}\leq \prob^{x,k-1} (\ret_k(0))
    \leq \frac{\E}{k} +\frac {C\log k}{k^2}.$$

\end{lemma}

Lemma \ref{returnafter} tells us that resampling a small
percentage of the moves is enough to get a very low correlation
between hitting $\dd$ before and after the re-randomization.
\begin{lemma} \label{returnafter}
There exists $C$ such that for all $k$, $n\geq
s_{k-1}+s_k/2^{10k}$, for all $I \subset \{1,\dots,n\}$ with
$|I|\geq s_k/2^{10k}$ and for all
 $\{x_i(t)\}_{i \in \{1,\dots,n\} \setminus I}$
$$    \prob \left(\exists j \in \{n,\dots,s_k\}
    \mbox{ such that } S_j(t)=\0\ \right |\
    \{x_i(t)\}_{i \in \{1,\dots,n\}\setminus I} )\leq
    \frac{C}{k}.$$
\end{lemma}

The last lemma we quote says there is a low correlation between
hitting $\dd$ (during some segment $M_k$) at different times.
\begin{lemma}  \label{bothreturn}
There exists $C$ such that for any $t$, any $k>\ko(t)$ and any
$x,y \in \ann_{k-1}$
$$\prob^{x,y,k-1} (\ret_k(0,t))
    \leq \frac{C}{k^2}.$$
\end{lemma}

\subsection{Estimating $\prob(SE_k(0))$ and $\prob(SE_k(0,t))$}

We start by giving estimates for the event $SR_j(t)$.
\begin{lem} \label{srj}
 There exists $C>0$, such that for any $j$ and any $x\in A_{j-1}$,

 $$\frac{\pi}{j^3} - \frac{C \log j}{j^4} \leq \prob^{x,j-1}(SR_j(0))
\leq \frac{\pi}{j^3} + \frac{C \log j}{j^4}.$$

%
%
%
\end{lem}
\begin{proof}
 Let $n$ be the first step in the segment $M_j$ such that
$S_{n}(0) \in \dd$, letting $n=\infty$ if no such step exists. By
Lemma \ref{mainlem},

\begin{equation}
\frac{2}{j} - \frac{C \log j}{j^2} \leq \prob^{x,j-1}(n<\infty) \leq
\frac{2}{j} + \frac{C \log j}{j^2}. \label{starstar}
\end{equation}

Let $B_i$ denote the event that the walk does not hit $\0$ between
step $i$ and $s_j$. Then
\begin{eqnarray}
\lefteqn{ \prob^{x,j-1}(SR_j(0)\cap G_j(0))}\hspace{.25in}&& \nonumber \\
&=& \sum_{s_{j-1} \leq i < s_j}\prob^{x,j-1}(n = i) \prob(B_i\cap
G_j(0) \,| \, n=i).\label {sss}
\end{eqnarray}
 For any $s_{j-1} \leq i < s_j$,
$B_i\cap G_j(0)$ is included in the event that after step $i$, the
walk reaches $A_j$ before reaching $\0$. Therefore by the second
statement in Lemma \ref{lawler} there exists $C$ such that for all
$i$
\begin{equation}
 \prob\left(B_i(0)\cap G_j(0) \,|\, n=i\right)\leq
\frac{\pi}{2j^2}  +\frac{C}{j^4}. \label{starstarstar}
\end{equation}
Putting (\ref{starstarstar}) and (\ref{starstar}) into (\ref{sss})
we get

\begin{eqnarray}
\lefteqn{\prob^{x,j-1}\left(SR_j(0)\cap G_j(0)\right)}\hspace{.75in}&&\nonumber\\
 &= &   \sum_i \prob^{x,j-1}(n=i)\cdot  \prob\left(B_i\cap G_j(0)\, | \, n=i\right)\nonumber\\
 &\leq &   \sum_i \prob^{x,j-1}(n=i)\cdot
        \sup_i\left(\prob\left(B_i\cap G_j(0)\, | \, n=i\right)\right)\nonumber\\
 &\leq &    \prob^{x,j-1}(n < \infty)\cdot\left(\frac{\pi}{2j^2}  +\frac{C}{j^4}\right)\nonumber\\
 & \leq &   \left(\frac{2}{j} + \frac{C \log j}{j^2}\right)\cdot
            \left(\frac{\pi}{2j^2} +  \frac{C}{j^4}\right)\nonumber\\
 &\leq&     \frac{\pi}{j^3} + \frac{C\log j}{j^4}. \label{sandor}
\end{eqnarray}
By Lemma \ref{gbound} we have that for all $x\in A_{j-1}$
 \begin{equation}
 \prob^{x,j-1}(G_j(0)^c) \leq \frac{C}{j^{10}}. \label {chuck}
 \end{equation}
Thus using (\ref{sandor}) and (\ref{chuck}) together with
conditional probabilities with respect to $G_j(0)$, we get

\begin{eqnarray*}
\prob^{x,j-1}(SR_j(0))
 &\leq&     \prob^{x,j-1}\left(SR_j(0)\cap G_j(0)\right) +
             \prob^{x,j-1}(G_j(0)^c)\\
 &\leq&     \frac{\pi}{j^3} + \frac{C\log j}{j^4}+\frac{C}{j^{10}}\\
 &\leq&     \frac{\pi}{j^3} + \frac{C\log j}{j^4}.
\end{eqnarray*}
which establishes the upper bound.

To get the lower bound, note that conditioned on the event $G_j(0)$,
if after step $i$ the walk reaches $A_j$ before reaching $(0,0)$,
and then does not hit $(0,0)$ in the next $s_j$ steps, then $B_i$
occurs.

 Letting $E_1$ denote the event that after step $i$ we reach
$A_j$ before reaching $(0,0)$, and letting $E_2$ denote the event
that a simple random walk, starting at $A_j$, does not hit $(0,0)$
in the next $s_j$ steps. We get, using the  Markovian property,
that
\begin{equation} \label{bn0lower}
\prob^{x,j-1}(B_i | \ G_j(0) \cap \{n = i\}) \geq \prob(E_1 \ | \
n=i) \prob(E_2).
\end{equation}

By Lemma \ref{lawler} we have
 \begin{equation}\label{E1}
 \prob(E_1 \ | \ n=i) \geq \left(\frac{\pi}{2j^2} -
 \frac{C}{j^4}\right).
 \end{equation}
 And by the Markovian property and stationarity,
since a walk starting at $A_j$ must hit $A_{j-1}$ before reaching
$0$, we can use Lemma \ref{mainlem} to bound
\begin{equation}\label{E2}
 \prob(E_2) \geq \left(1-\frac{2}{j} - \frac{C\log j}{j^2}\right).
\end{equation}

Combining (\ref{E1}) and  (\ref{E2}) into (\ref{bn0lower}) we get
\begin{eqnarray}
 \prob^{x,j-1}(B_i \ | \ G_j(0) \cap \{n = i\} )
 &\geq & \left(\frac{\pi}{2j^2} - \frac{C}{j^4}\right)\left(1-\frac{2}{j}
        - \frac{C\log j}{j^2}\right) \nonumber \\
 &\geq&  \frac{\pi}{2j^2} - \frac{C}{j^3}. \label{dumb}
\end{eqnarray}

Using (\ref{dumb}) for all $j$ and $x \in A_{j-1}$
\begin{eqnarray}
\prob^{x,j-1}(SR_j(0) \, | \, G_j(0))
& =    & \sum_i\prob^{x,j-1}(n =i)
    \prob(B_i \, | \, G_j(0)\cap\{n=i\}) \nonumber \\
& \geq & \prob^{x,j-1}(n < \infty)\inf_i
    \prob(B_i\, |  G_j(0)\cap\{n=i\}) \nonumber \\
& \geq & \left(\frac{2}{j}-\frac{C\log j}{j^2}\right)\cdot
        \left(\frac{\pi}{2j^2} - \frac{C}{j^3}\right)\nonumber \\
&\geq & \frac{\pi}{j^3} - \frac{C\log j}{j^4}. \label{line34}
\end{eqnarray}

Thus for all $j$ and $x \in A_{j-1}$, conditioning on $G_j(0)$ and
using Lemma \ref{gbound} and (\ref{line34}) implies that

\begin{eqnarray*} \prob^{x,j-1}(SR_j(0))
 &\geq& \prob^{x,j-1}(SR_j(0) \,| \, G_j(0)) - \prob^{x,j-1}(G_j(0)^c) \\
 &\geq& \frac{\pi}{j^3} - \frac{C\log j}{j^4} - \frac{C}{j^{10}} \\
 &\geq& \frac{\pi}{j^3} - \frac{C\log j}{j^4}.
 \end{eqnarray*}

\end{proof}

The next two lemmas give a lower bound on $\prob(SE_k(0))$. We
define \be \beta_k := \prod_{2^k \leq j <
2^{k+1}}\left(1-\frac{2}{j}\right).
 \label{betak} \ee
 This is an approximation of the probability of never hitting $\dd$
 during the super segment $W_k$.
 Note that there exists $c>0$ such that $c < \beta_k <1$ for all $k>1$.

\begin{lemma} \label{sekj}
For any $k$ and any $2^k \leq j < 2^{k+1}$, and any $x\in
A_{2^k-1}$
$$\prob^{x,2^k-1}\left(SE_k^j(0)\right) \geq \beta_k \frac{\pi}{j^3} - \frac{C \log
j}{j^4}.$$

\end{lemma}
\begin{proof}

By the Markovian property, for any $x\in A_{2^k-1}$
 \pagebreak[3]

\begin{eqnarray}\label{sekmarkov}
\lefteqn{\prob^{x,2^k-1}\Big(SE_k^j(0) \Big)}  && \nonumber \\
 &\geq& \min_{x\in A_{j-1}}\prob^{x,j-1}\Big(SR_j(0) \cap G_j(0)\Big) \cdot \nonumber\\
 && \cdot \prod_{2^k \leq l<2^{k+1}, l\neq j}
        \left(\min_{x\in A_{l-1}} \prob^{x,l-1}\left((R_l(0)\cap G_l(0)\right)^c\right)
        \nonumber  \\
 & \geq & \min_{x\in A_{j-1}}\prob^{x,j-1}\Big(SR_j(0)\Big) \prod_{2^k \leq l<2^{k+1}, l\neq j}
        \left(\min_{x\in A_{l-1}}
        \prob^{x,l-1}\left(R_l(0)\right)^c\right) \nonumber \\
        && - \sum_{2^k \leq l <
        2^{k+1}} \prob^{x,l-1}\Big((G_l)^c\Big)
\end{eqnarray}

By Lemma \ref{srj}
 \begin{equation}\label{lsek1}
 \min_{x\in
A_{j-1}}\prob^{x,j-1}(SR_j(0)) \geq \frac{\pi}{j^3} -\frac{C \log
j}{j^4}.
 \end{equation}
 By Lemma \ref{mainlem}
 \[\min_{x\in A_{l-1}}\prob^{x,l-1}((R_l(0))^c) \geq
\left(1-\frac{2}{l}\right)-\frac{C \log l}{l^2}.\] Multiplying
over all $l$ such that $2^k \leq l<2^{k+1}$ and  $l\neq j$, we get
\begin{eqnarray}
\lefteqn{\prod_{2^k \leq l<2^{k+1}, l\neq j}
            \left(\min_{x\in A_{l-1}}
    \prob^{x,l-1}((R_l(0))^c)\right)} \hspace{.5in}&& \nonumber\\
 & \geq&  \prod_{2^k \leq l<2^{k+1}, l\neq j}\left(1-\frac{2}{l} -
    \frac{C \log l}{l^2}\right) \nonumber \\
 & \geq& \prod_{2^k \leq l<2^{k+1}, l\neq j}\left(1-\frac{2}{l}\right)
 \nonumber \\
    && \ \  - \,  \sum_{2^k \leq l<2^{k+1},l\neq j} \left(\frac{C \log l}{l^2}\prod_{2^k \leq m<2^{k+1}, m\neq
    l,j}\left(1-\frac{2}{m}\right)\right) \nonumber \\
    & \geq& \prod_{2^k \leq l<2^{k+1}}\left(1-\frac{2}{l}\right) \ - \
    \sum_{2^k \leq l<2^{k+1}} \frac{C \log l}{l^2}\prod_{2^k \leq l<2^{k+1}}\left(1-\frac{2}{l}\right) \nonumber \\
 & \geq& \prod_{2^k \leq l<2^{k+1}}\left(1-\frac{2}{l}\right) \ - \
    2^k \frac{C\log 2^k}{2^{2k}}\prod_{2^k \leq l<2^{k+1}}\left(1-\frac{2}{l}\right) \nonumber \\
&\geq& \beta_k - \frac{\beta_k C \log 2^k}{2^k} \nonumber\\
&\geq& \beta_k - \frac{C\log j}{j}. \label{lsek2}
 \end{eqnarray}

Last, by Lemma \ref{gbound}
\begin{equation}\label{lsek3}
\sum_{2^k \leq l < 2^{k+1}} \prob^{x,l-1}((G_l)^c)  \leq 2^k
\frac{C}{2^{10k}} \leq \frac{C}{2^{9k}}.
 \end{equation}

 Putting (\ref{lsek1}),(\ref{lsek2}) and (\ref{lsek3}) into (\ref{sekmarkov}) we
get
\begin{eqnarray*}
\lefteqn{\prob^{x,2^k-1}\bigg(SE_k^j(0) \bigg) } \hspace{.5in}&& \\
 &\geq& \min_{x\in A_{j-1}}\prob^{x,j-1}\Big(SR_j(0)\Big) \prod_{2^k \leq l<2^{k+1}, l\neq j}
        \left(\min_{x\in A_{l-1}}
        \prob^{x,l-1}\left(R_l(0)\right)^c\right)   \\
        &&  - \sum_{2^k \leq l < 2^{k+1}} \prob^{x,l-1}\Big((G_l)^c\Big)
        \nonumber\\
 &\geq&
  \left( \frac{\pi}{j^3} -\frac{C \log j}{j^4} \right)
  \left( \beta_k - \frac{C \log j}{j} \right) - \frac{C}{2^{9k}}\nonumber \\
  &\geq&  \beta_k \frac{\pi}{j^3} - \frac{C \log j}{j^4}.
\end{eqnarray*}

\end{proof}


\begin{lemma} \label{sek}
\[\prob\bigg(SE_k(0) \ | \ SE_{k-1}(0)\bigg)\geq \pi \beta_k \left(\sum_{2^k \leq j <
2^{k+1}}\frac{1}{j^3}\right) \   - \frac{C k}{2^{3k}}.\]
\end{lemma}
\begin{proof}
$SE_{k-1}(0) \subseteq G_{2^k-1}(0)$, and therefore by the Markovian
property of simple random walk, we have
$$\prob(SE_k(0) \ | \ SE_{k-1}(0)) \geq \min_{x\in A^{2^k-1}}\prob^{x,2^k-1}(SE_k(0)).$$

By the definition of $SE_k^j$, the events $SE_k^{j_1}$ and
$SE_k^{j_2}$ are disjoint  for $j_1\neq j_2$, so \pagebreak[1]
\begin{eqnarray} \min_{x\in
A^{2^k-1}}\prob^{x,2^k-1}(SE_k(0))
&=&  \min_{x\in A^{2^k-1}}\sum_{2^k \leq j < 2^{k+1}}\prob^{x,2^k-1}(SE_k^j(0)) \nonumber \\
& \geq & \ \sum_{2^k \leq j < 2^{k+1}}\left(\min_{x\in
A^{2^k-1}}\prob^{x,2^k-1}(SE_k^j(0))\right) \nonumber \\
& \geq & \ \sum_{2^k \leq j < 2^{k+1}}\Big(\beta_k \frac{\pi}{j^3}
- \frac{C \log j}{j^4}\Big) \label{lock} \\
& \geq & \ \pi \beta_k \left(\sum_{2^k \leq j <
2^{k+1}}\frac{1}{j^3} \right) \ - \  \frac{C k}{2^{3k}}. \nonumber
 \end{eqnarray}
 The inequality in (\ref{lock}) follows from Lemma
 \ref{sekj}.
\end{proof}

$\ $\\ Next we will work to bound from above the probability of
$SE_k(0,t)$.

First, we deal with the case that the walk hits $\dd$ in different
sub-segments of $W_k$ in times $0$ and $t$.
\begin{lemma} \label{sek12}
There exists a constant $C>0$ such that for any $k> K'(t)$ and any
$j_1\neq j_2$, $2^k \leq j_1,j_2 < 2^{k+1}$, we have $$\max_{x,y\in
A_{2^k-1}}\prob^{x,y,2^k-1}(SE_k^{j_1,j_2}(0,t)) \leq \beta_k^2
\frac{\pi}{j_1^3} \frac{\pi}{j_2^3} + \frac{C k}{2^{7k}}.$$
\end{lemma}

\begin{proof}
Recall that

\begin{eqnarray*}
\lefteqn{SE_k^{j_1,j_2}(0,t) = }& & \\ & & \bigg\{SR_{j_1}(0) \cap
(R_{j_1}(t))^c\bigg\}\cap\bigg\{(R_{j_2}(0))^c \cap
SR_{j_2}(t)\bigg\} \cap \\ & & \quad \bigg\{\bigcap_{2^k\leq l <
2^{k+1}, l\neq j_1,j_2}\Big((R_l(0))^c \cap (R_l(t))^c \,
\Big)\bigg\} \cap \bigg\{
 \bigcap_{2^k\leq l < 2^{k+1}}(G_l(0,t))^c\bigg\}.
\end{eqnarray*}

So by the Markovian property,
\begin{eqnarray}\label{easy1}
\lefteqn{\max_{x,y\in A_{2^k-1}}\prob^{x,y,2^k-1}(SE_k^{j_1,j_2}(0,t)) }  \hspace{.5in}&& \nonumber \\
 & \leq & \max_{x,y\in A_{j_1-1}}\prob^{x,y,j_1-1}\Big(SR_{j_1}(0)\cap
(R_{j_1}(t))^c\Big)  \nonumber \\
& & \cdot \max_{x,y\in
    A_{j_2-1}}\prob^{x,y,j_2-1}\bigg((R_{j_2}(0))^c
    \cap SR_{j_2}(t)\bigg)  \\
& & \cdot \prod_{2^k\leq l < 2^{k+1}, l\neq j_1,j_2}\max_{x,y\in
A_{l-1}}\prob^{x,y,l-1}\bigg((R_l(0))^c \cap (R_l(t))^c\bigg) \nonumber \\
&& \ + \ \max_{x,y \in A_{l-1}}
 \prob^{x,y,2^k-1}\left(\left(\bigcap_{2^k\leq l <
2^{k+1}}G_l(0,t)\right)^c\right). \nonumber
\end{eqnarray}

Now we will estimate each of the four terms in (\ref{easy1}).

\begin{eqnarray}
\max_{x,y\in A_{j_1-1}}
    \prob^{x,y,j_1-1}(SR_{j_1}(0)\cap(R_{j_1}(t))^c)
 &\leq &  \max_{x\in A_{j_1-1}}\prob^{x,j_1-1}(SR_{j_1}(0)) \nonumber\\
 &\leq & \frac{\pi}{j_1^3} + \frac{C \log j_1}{j_1^4}. \label{easy2}
\end{eqnarray}
The last inequality holding by Lemma \ref{srj}.

Similarly, by symmetry between times $0$ and $t$,
\begin{equation}\label{easy3}
\max_{x,y\in A_{j_2-1}}\prob^{x,y,j_2-1}\Big((R_{j_2}(0))^c \cap
SR_{j_2}(t)\Big)\leq \frac{\pi}{j_2^3} + \frac{C \log j_2}{j_2^4}.
\end{equation}

To estimate the third term in (\ref{easy1})  we  use Lemma
\ref{bothreturn} to estimate each term in the product.  Lemma
\ref{bothreturn} says that there exists $C>0$ such that for any
$l$ and any $x,y\in A_{l-1}$
$$\prob^{x,y,l-1}(R_l(0,t)) \leq \frac{C}{l^2}.$$
And the lower bound from Lemma \ref{mainlem}:
$$\prob^{x,l-1}(R_l(0))\geq \frac{2}{l} - \frac{C\log l}{l^2}.$$

 We then get that
\begin{eqnarray*}
\lefteqn{\max_{x,y\in A_{l-1}} \prob^{x,y,l-1}\bigg((R_l(0))^c \cap (R_l(t))^c\,\bigg) }\hspace{.5in} &&\\
 &\leq & 1 - \min_{x,y\in A_{l-1}}\prob^{x,y,l-1}(R_l(0)) \\
        && \ \  - \min_{x,y\in A_{l-1}}\prob^{x,y,l-1}(R_l(t)) +
        \max_{x,y,l-1}\prob^{x,y,l-1}(R_l(0,t)) \\
        &\leq & 1 - 2\min_{x\in A_{l-1}}\prob^{x,l-1}(R_l(0)) +
        \max_{x,y,l-1}\prob^{x,y,l-1}(R_l(0,t)) \\
 & \leq & 1 - \frac{4}{l} +
\frac{C \log l}{l^2}.
\end{eqnarray*}

 Which yields
\begin{eqnarray}\label{easy4}
\lefteqn{\prod_{2^k \leq l < 2^{k+1}, l\neq j_1,j_2}
\max_{x,y,l-1} \prob^{x,y,l-1}\Big((R_l(0))^c \cap (R_l(t))^c\Big)}  \hspace{.5in}&& \nonumber \\
 & \leq &  \prod_{2^k \leq l < 2^{k+1}, l\neq
    j_1,j_2}\Big(1-\frac{4}{l}+\frac{C\log l}{l^2}\Big)  \nonumber \\
 & \leq &  \prod_{2^k \leq l < 2^{k+1}, l\neq
    j_1,j_2}\Big(1-\frac{4}{l}+\frac{4}{l^2} + \frac{C\log l}{l^2}\Big)  \nonumber \\
 &\leq &  \prod_{2^k \leq l < 2^{k+1}, l\neq j_1,j_2}
    \Big(1-\frac{4}{l}+\frac{4}{l^2}\Big) \nonumber\\
 && + \ \sum_{2^k \leq l < 2^{k+1}, l\neq j_1,j_2}\left(\frac{C \log l}{l^2}
 \prod_{2^k \leq m < 2^{k+1},m\neq l,j_1,j_2}\Big(1-\frac{4}{m}+\frac{4}{m^2}\Big)\right) \nonumber \\
& \leq & \frac{\beta_k^2}{(1-\frac{2}{j_1})^2 (1-\frac{2}{j_2})^2}
+ \sum_{2^k \leq l < 2^{k+1}, l\neq j_1,j_2}\frac{C \log l}{l^2}
\nonumber \\
& \leq & \beta_k^2 +\frac{C k}{2^k}.
\end{eqnarray}

For the fourth factor in (\ref{easy1}), Corollary \ref{gbigbound}
gives
\begin{eqnarray}\label{easy5}
\max_{x,y \in A_{l-1}}\lefteqn{\prob^{x,y,2^k-1}
 \left(\bigcap_{2^k \leq l < 2^{k+1}}G_l(0,t)\right)^c\, }\hspace{.5in} && \nonumber \\
 &\leq & \prob^{x,2^k-1}\left(\bigcap_{2^k \leq l < 2^{k+1}}G_l(0)\right)^c\ +
 \prob^{y,2^k-1}\left(\bigcap_{2^k \leq l <
 2^{k+1}}G_l(0)\right)^c\,\nonumber \\
 &\leq& \frac{C}{2^{9k}}.
\end{eqnarray}

Combining (\ref{easy2}), (\ref{easy3}), (\ref{easy4}) and
(\ref{easy5}) into (\ref{easy1}) we get
\begin{eqnarray*}
\lefteqn{\max_{x,y\in
A_{2^k-1}}\prob^{x,y,2^k-1}(SE_k^{j_1,j_2}(0,t))} \hspace{.5in} &&\\
& \leq &
     \bigg(\frac{\pi}{j_1^3} + \frac{C \log j_1}{j_1^4}\bigg)\cdot
     \bigg(\frac{\pi}{j_2^3} + \frac{C \log j_2}{j_2^4}\bigg) \cdot
     \bigg(\beta_k^2 + \frac{C k}{2^k}\bigg) + \frac{C}{2^{9k}} \\
& \leq &
     \bigg(\frac{\pi}{j_1^3} + \frac{C k}{2^{4k}}\bigg)\cdot
     \bigg(\frac{\pi}{j_2^3} + \frac{C k}{2^{4k}}\bigg) \cdot
     \bigg(\beta_k^2 + \frac{C k}{2^k}\bigg) + \frac{C}{2^{9k}} \\
& \leq &  \beta_k^2 \frac{\pi}{j_1^3} \frac{\pi}{j_2^3} + \frac{C
k}{2^{7k}}.
\end{eqnarray*}
\end{proof}


Next, we deal with the case in which the walk hits $\dd$ in the same
sub-segment $M_j\subset W_k$, and show it holds with only negligible
probability.
\begin{lem}\label{hard_case}
There exists $C>0$ such that for any $t>0$, $k>K'(t)$, \\ $2^k \leq
j < 2^{k+1}$ and $x,y\in A_{j-1}$
\[P^{x,y,j-1}\big(SR_j(0,t)\big) \leq \frac{C}{j^6}.\]

\end{lem}
\begin{proof}
 Let $n_0$ be the first step $\dd$ is hit in the segment  $M_j$
 at
time $0$, and $n_t$ be the first step that $\dd$ is hit in $M_j$
at time $t$, letting them equal $\infty$ if $\dd$ is not hit in
$M_j$ at all at that time. Notice that by symmetry, we can
estimate
$$\prob^{x,y,j-1}\big(SR_j(0,t)\big) \leq 2\prob^{x,y,j-1}\big(SR_j(0,t) \cap (n_0 \leq n_t)\big).$$

Let $$r_1= \frac{2^{(j-1)^2}}{2^{2(j+2)}j^{12}}$$ and
 $$\tau_{r_1} = \min\{ m>n_0  \ : \ S_m(0) > r_1\},$$

Denote by $L_1$ the event that $G_j(0)$ occurs, $\tau_{r_1} <
\infty$ and
 $S_m(0)\neq (0,0)$ for all integers $m$ between $n_0$ and $\tau_{r_1}$. \\
Denote by  $L_2$ the event that $G_j(t)$ occurs and $S_m(t)\neq
 (0,0)$ for all integers $m$ between $n_t$ and $s_j$.\\

Then
 $$SR_j(0) \cap G_j(0)\subseteq R_j(0)\cap L_1$$
and
 $$SR_j(t)\cap G_j(t) \cap \{n_0 \leq n_t\} \subseteq
R_j(t)\cap L_2.$$ Combining these we get
\begin{eqnarray*}
\lefteqn{ SR_j(0,t) \cap G_j(0,t)\cap \{n_0\leq n_t \}}\hspace{1in}&&\\
 &\subset & \bigg( SR_j(0)\cap G_j(0) \bigg) \cap \bigg(SR_j(t) \cap G_j(t) \cap \{n_0\leq n_t \}\bigg)\\
 &\subset & R_j(0)\cap L_1 \cap R_j(t) \cap L_2 .
\end{eqnarray*}
Thus we get

\begin{eqnarray}\label{hardestimate}
\lefteqn{\prob^{x,y,j-1}\big(SR_j(0,t)\big)} \nonumber && \\
 &\leq  & 2\prob^{x,y,j-1}\big(SR_j(0,t)\cap \{n_0\leq n_t\}\big) \nonumber \\
 &\leq  & 2\prob^{x,y,j-1}\big(SR_j(0,t)\cap G_j(0,t)
    \cap \{n_0 \leq n_t\}\big)  + 2\prob^{x,y,j-1}((G_j(0,t))^c)\nonumber  \\
 & \leq & 2\prob^{x,y,j-1}\big(R_j(0)\big) \cdot
\prob^{x,y,j-1}\big(L_1 \ | \ R_j(0)\big) \cdot
\prob^{x,j-1}\big(R_j(t)\ | \ R_j(0)\cap L_1\big) \cdot \nonumber \\
 && \quad
  \cdot \prob^{x,y,j-1}\big(L_2 \ | \ R_j(0,t) \cap L_1 \big)
  + \ 2\prob^{x,y,j-1}((G_j(0,t))^c).
\end{eqnarray}

To bound the probabilities of the five terms on the right hand side
of (\ref{hardestimate}) we now introduce six new events, bound the
probabilities of these events and then use these bounds to bound the
terms in (\ref{hardestimate}).

Let $B_1$  be the event that
$$n_0 < s_{j-1} +
\frac{2^{2(j-1)^2}}{j^{6}}.$$ Recall that by Corollary \ref{hit0},
for any $M,N>0$, if $|S_{s_{j-1}(0)}|>M$, then
\[
\prob(\exists s_{j-1}<n'<s_{j-1}+N \ : \  |S_{n'}(0)|\leq 1) \leq
\frac{CN}{M^2}.\]

 Since for any $x\in A_{j-1}$
we have $|x| \geq 2^{(j-1)^2}$, we can set $M=2^{(j-1)^2}$ and
$N=\frac{2^{2(j-1)^2}}{j^{6}}$ to conclude

\begin{eqnarray}\label{b1bound}
\prob^{x,y,j-1}(B_1) & = &\prob^{x,y,j-1} \left(\exists
s_{j-1}<n'<s_{j-1}+ \frac{2^{2(j-1)^2}}{j^{12}} : \
|S_{n'}(0)|\leq 1 \right)    \nonumber
\\  & \leq & \frac{C\frac{2^{2(j-1)^2}}{j^{6}}}{(2^{(j-1)^2})^2} \nonumber \\
&\leq& \frac{C}{j^6}.
\end{eqnarray}

Let $I$ denote the set of all indices between $s_{j-1}$ and $
s_{j-1}+\frac{2^{2(j-1)^2}}{j^{6}}$ for which, conditioned on our
Poisson process, $X_i(0)$ and $X_i(t)$ are independent.\\
Let $B_2$ be the event that $|I| < \frac{2^{2(j-1)^2}}{2^{j+2}j^6}$.\\

If $B_1$ does not occur, then $|I|$ is the sum of at least
$\frac{2^{2(j-1)^2}}{j^{6}}$ i.i.d. indicator variables each
happening with probability
 $$p = 1-e^t> \frac{1}{2}\min(1,t) \geq \frac{1}{2^{j+1}},$$
(the last inequality holds because $j \geq 2^k \geq K(t) \geq
|\log t|$). Therefore by monotonicity in $t$ and in the number of
indicators, $|I|$ stochastically dominates the sum of
$N=\frac{2^{2(j-1)^2}}{j^{6}}$ i.i.d. indicator variables $T_i$
each equaling $1$ with probability $p=\frac{1}{2^{j+1}}$. Let
$T=\sum_{1\leq i \leq N}T_i$.

As we are conditioning on $(B_1)^c$ we have that $|I|$  dominates
$T$ and if $T\geq \frac{NP}{2}$ then

$$|I|\geq \frac{NP}{2}=\frac{2^{2(j-1)^2}}{2^{j+2}j^6}$$
and $B_2$ does not occur.

 By the Chernoff bound
$$\prob\bigg(T<\frac{Np}{2}\bigg) \leq \prob\bigg(|T-Np|\geq \frac{Np}{2}\bigg) \leq
2e^{-cNp}$$ for some absolute constant $c>0$.  Thus

\begin{equation} \label{cor7}
\prob^{x,y,j-1}\big(B_2 \ | \ (B_1)^c\big) \leq \prob
\bigg(T>\frac{Np}{2}\bigg)
 \leq 2e^{-\frac{c2^{2(j-1)^2}}{2^{j+2} j^{6}}} \leq \frac{C}{j^6}.
\end{equation}

Put $r =\frac{2^{(j-1)^2}}{2^{j+2}j^6}$.

Let $B_3$ be the event that $|S_{n_0}(t)| < r$ and
  $$S_{n_0}(t) =\sum_{i<n_0}X_i(t).$$
  Rearranging the order of summation we can
write
 $$S_{n_0}(t) = \sum_{i<n_0}X_i(t) = \sum_{i<n_0, i\notin I}X_i(t) + \sum_{i\in I}X_i(t).$$
  Put
   $$x = \sum_{i<n_0, i\notin I}X_i(t).$$
Since all variables $\{X_i\}_i\in I$ have been re-rolled between
time $0$ and $t$, the probability that $|S_{n_0}(t) < r|$ is the
same as the probability that simple random walk starting at $x$
will be at distance less than $r$ from $(0,0)$ after $|I|$ steps.
Therefore by the second part of Lemma \ref{leave}
\begin{equation} \prob^{x,y,j-1}(S_{n_0}(t) <
\frac{\sqrt{|I|}}{j^3} ) \leq \frac{C}{j^6}
\end{equation}
And since $(B_2)^c$ implies $|I|> \frac{2^{2(j-1)^2}}{2^{j+2}j^6}
 > r^2 j^6$ we get
\begin{equation} \prob^{x,y,j-1}(B_3 \ | (B_2)^c) < \frac{C}{j^6}
\end{equation}

Next we bound the probability that $n_t-n_0$ is small. Let $B_4$
be the event that $0< n_t-n_0 < \frac{r^2}{j^6}
 = \frac{2^{2(j-1)^2}}{2^{j+2}j^12} $. Conditioning on
$|S_{n_0}(t)|> r$, $B_4$ implies that the random walk at time $t$
gets to distance $r$ from its position at step $n_0$ in less than
$\frac{r^2}{j^6}$ steps, therefore we can apply the first part of
Lemma \ref{leave} to bound
\begin{equation}
\prob^{x,y,j-1}( B_4 | \ (B_3)^c) < \frac{C}{j^6}.
\end{equation}

Let $I_1$ denote the set of all indices  between $n_0$ and $n_t$
for which, conditioned on our Poisson process, $X_i(0)$ and
$X_i(t)$ are independent. Let $B_5$ denote the event that $|I_1| <
\frac{2^{2(j-1)^2}}{2^{4(j+2)}j^12}$. A similar calculation as
done for $B_2$ shows that $\prob^{x,y,j-1}(B_5) < \frac{C}{j^6}$

Let $B_6$ be the event that $|S_{n_t}(0)| < r_1$. Using Lemma
\ref{leave} in a calculation similar to the one for $B_3$ gives
that
\begin{equation} \label{nopunc}
\prob^{x,y,j-1}(B_6 \ | (B_5)^c) < \frac{C}{j^6}.
\end{equation}

 Thus we get
\begin{equation}
\prob^{x,y,j-1}(B_i) \leq \frac{C}{j^6} \quad\quad i=1,\dots,6.
\end{equation}

We will estimate the five probabilities in (\ref{hardestimate}) as follows: \\
\begin{enumerate}
\item By Lemma \ref{mainlem},
\begin{equation}\label{rj0} \prob^{x,y,j-1}\big(R_j(0)\big)=\prob^{x,j-1}\big(R_j(0)\big) \leq \frac{C}{j}.
\end{equation}

\item $L_1$ is included in the event that after step $n_0$, at for
which $|S_{n_0}(0)|=1$, the walk reaches distance $r$ before
hitting $(0,0)$, therefore by Lemma \ref{lawler},
\begin{equation}\label{L1}\prob^{x,y,j-1}\big(L_1 \ | \ R_j(0)\big)
\leq \frac{C}{\log r_1} \leq  \frac{C}{j^2}.
\end{equation}

\item  To estimate $\prob^{x,y,j-1}\big(R_j(t) \ | \ R_j(0) \cap
L_1\big)$, notice that $R_j(0)$ and $L_1$ depend only on what
happens at time $0$, and $(B_2)^c$ implies $|I|\geq
\frac{s_j}{2^{10j}}$ (For all large enough $j$) , so using Lemma
\ref{returnafter} we get
\[\prob^{x,y,j-1}\big(R_j(t) \ |\ (B^2)^c \cap R_j(0) \cap L_1\big) \leq
 \frac{C}{j},\]
 and therefore
\begin{equation}\label{rjt}\prob^{x,y,j-1}\big(R_j(t) \ | \ R_j(0) \cap L_1\big) \leq
\frac{C}{j} + \frac{C}{j^6} \leq \frac{C}{j}.\end{equation}

\item To estimate $\prob^{x,y,j-1}\big(L_2 \ | \ R_j(0,t)\cap
L_1\big) $, we partition according to the value of $n_t$. Since
for a given $n_t$, the event $L_2$ is independent of the variables
$\{X_i(0)\}_{i<n_t} \cup \{X_i(t)\}_{i<n_t}$, and the event
$R_j(0,t)\cap L_1 \cap (B_6)^c$ depend only on these variables, we
have $\forall N_t$
\[ \prob^{x,y,j-1}\big(L_2 \ | \ R_j(0,t)\cap L_1 \cap (B_6)^c \cap (n_t=N_t)\big) =
\prob^{x,y,j-1}\big(L_2 \ | \ n_t=N_t\big).\]

Since by Lemma \ref{lawler}, $$\prob^{x,y,j-1}(L_2 \ | \ n_t=N_t)
\leq \frac{C}{j^2},$$ for any value of $N_t$, we deduce that
$$ \prob^{x,y,j-1}\big(L_2 \ |\ \ R_j(0,t)\cap L_1 \cap (B_6)^c \big)
\leq \frac{C}{j^2}$$
 and therefore
\begin{equation}\label{L2}\prob^{x,y,j-1}\big(L_2 \ |\ \ R_j(0,t)\cap L_1 \big) \leq
\frac{C}{j^2} + \prob^{x,y,j-1}(B_6) \leq \frac{C}{j^2}
.\end{equation}

\item Last we use Lemma \ref{gbound} to bound
\begin{equation} \label{gj} \prob^{x,y,j-1}((G_j(0,t))^c) \leq 2\prob^{x,j-1}((G_j(0))^c) \leq \frac{C}{j^{10}}. \end{equation}
\end{enumerate}

Combining (\ref{rj0}), (\ref{L1}), (\ref{rjt}), (\ref{L2}) and
(\ref{gj}) into (\ref{hardestimate}) we get

\begin{eqnarray*}
\lefteqn{\prob^{x,y,j-1}\big(SR_j(0,t) \big) }&&\\
  &\leq & 2\prob^{x,y,j-1}\big(R_j(0)\big) \cdot \prob^{x,y,j-1}\big(L_1 \ | \ R_j(0)\big) \\
        && \ \cdot \, \prob^{x,y,j-1}\big(R_j(t)\
            | \ R_j(0)\cap L_1\big) \cdot   \prob^{x,y,j-1}\big(L_2 \ | \ R_j(0,t) \cap L_1\big) \\
            &&  + \ \prob^{x,y,j-1}((G_j(0,t))^c)\\
  &\leq &  2\frac{C}{j} \cdot \frac{C}{j^2} \cdot \frac{C}{j} \cdot \frac{C}{j^2} + \frac{C}{j^{10}}\\
  & \leq & \frac{C}{j^6}.
\end{eqnarray*}


\end{proof}

\begin{lemma}\label{sek0t}
There exists a constant $C>0$ s.t. for any $t>0$ and any $k>K'(t)$
$$\prob(SE_k(0,t) \ | \ SE_{k-1}(0,t)) \leq \left(\beta_k \pi \sum_{2^k \leq j
< 2^{k+1}} \frac{1}{j^3}\right)^2 + \frac{Ck}{2^{5k}}. $$
\end{lemma}
\begin{proof}
$SE_{k-1}(0,t) \subseteq G_{2^k-1}(0,t)$, so by the Markovian
properties of simple random walk we have $$\prob(SE_k(0,t) \ | \
SE_{k-1}(0,t)) \leq \max_{x,y\in
A_{2^k-1}}\prob^{x,y,2^k-1}(SE_k(0,t)).$$
 $SE_k(0,t)$ is the disjoint union of the events $SE_k^{j_1,j_2}(0,t)$, therefore

\begin{eqnarray} \lefteqn{ \max_{x,y\in
A_{2^k-1}}\prob^{x,y,2^k-1}(SE_k(0,t))} \hspace{1in}&& \nonumber \\
&\leq & \sum_{2^k \leq j_1,j_2 < 2^{k+1}}\max_{x,y\in
        A_{2^k-1}}\prob^{x,y,2^k-1}(SE_k^{j_1,j_2}(0,t)) \nonumber \\
&=&\sum_{2^k \leq j_1\neq j_2 < 2^{k+1}}\max_{x,y\in
        A_{2^k-1}}\prob^{x,y,2^k-1}(SE_k^{j_1,j_2}(0,t)) \nonumber \\
&& + \ \sum_{2^k \leq j \leq 2^{k+1}} \max_{x,y\in
        A_{2^k-1}}\prob^{x,y,2^k-1}(SE_k^{j,j}(0,t)).
        \label{berry3}
\end{eqnarray}

By Lemma \ref{sek12}, for any $2^k \leq j_1\neq j_2 < 2^{k+1} $
\begin{equation} \label{berry1}
\max_{x,y\in
        A_{2^k-1}}\prob^{x,y,2^k-1}(SE_k^{j_1,j_2}(0,t)) \leq \beta_k^2
\frac{\pi}{j_1^3} \frac{\pi}{j_2^3} + \frac{C k}{2^{7k}}.
\end{equation}
$SE_k^{j,j}(0,t) \subseteq SR_j(0,t) $, therefore by Lemma
\ref{hard_case}, for any $2^k \leq j < 2^{k+1}$

\begin{equation} \label{berry2} \max_{x,y\in
        A_{2^k-1}}\prob^{x,y,2^k-1}(SE_k^{j,j}(0,t)) \leq
        \frac{C}{j^6}. \end{equation}

 Combining (\ref{berry1}) and (\ref{berry2}) into (\ref{berry3}) we conclude

\begin{eqnarray*} \lefteqn{ \max_{x,y\in
A_{2^k-1}}\prob^{x,y,2^k-1}(SE_k(0,t))}\hspace{.5in}  && \\
&& \leq  \ \beta_k^2\pi^2\sum_{2^k \leq j_1\neq j_2 < 2^{k+1}}
\left(\frac{1}{j_1^3}\frac{1}{j_2^3} + \frac{C k}{2^{7k}}\right) +
\sum_{2^k
\leq j < 2^{k+1}}\frac{C}{j^6} \\
&& \leq \ \left(\beta_k \pi \sum_{2^k \leq j < 2^{k+1}}
\frac{1}{j^3}\right)^2 + \frac{Ck}{2^{5k}}.
\end{eqnarray*}

\end{proof}

\subsection{Proof of Theorem \ref{fmt}}

\begin{pfof}{Theorem \ref{fmt}}
 Define
$$f(t,M) =\frac{\prob(\bigcap_{1\leq k \leq M}SE_k(0,t))}
            {(\prob\bigcap_{1\leq k \leq M}(SE_k(0)))^2}.$$

Set  
 $k_0 = K'(t)$. Then
\begin{equation}\label{eqfm}
f(t,M) \leq \frac{1}{\prob(SE_{k_0}(0))^2} \prod_{k=k_0+1}^M
\frac{\prob(SE_k(0,t) \ | \ SE_{k-1}(0,t))}{\prob(SE_k(0) \ | \
SE_{k-1}(0))^2}.
\end{equation}

We start by estimating the first factor:
\[ \frac{1}{\prob(SE_{k_0}(0))^2} \leq C\prod_{0 <  k \leq
 k_0}\frac{1}{\prob(SE_k(0) \ | \ SE_{k-1}(0))^2}. \]
 Using the lower bound from Lemma \ref{sek} we get that this is bounded above
by
\begin{eqnarray*}
 \frac{1}{\prob(SE_{k_0}(0))^2}
    &\leq & C\prod_{0 < k \leq k_0}\frac{1}{\left(\left(\beta_k \pi
        \sum_{2^k \leq j < 2^{k+1}}\frac{1}{j^3}\right) - \frac{C k}{2^{3k}}\right)^2}\\
 &\leq & C\prod_{0 < k \leq k_0} \frac{1}{(\frac{C}{2^{2k}} -
 \frac{C k}{2^{3k}})^2}   \\
    &\leq & C\prod_{0 < k \leq k_0}C2^{4k}  \\
    &\leq & C^{k_0} 2^{2k_0^2}  \\
    &\leq & C 2^{3k_0^2}.  \\
\end{eqnarray*}

We now recall that we chose $k_0$ to be the smallest integer larger
then $\log (|\log t| +1) + 1$, to get
   \begin{eqnarray}
 \frac{1}{\prob(SE_{k_0}(0))^2}
 &\leq&  C 2^{3k_0^2}\nonumber\\
 &\leq&  C \left(2^{k_0}\right)^{3k_0}\nonumber\\
 &\leq&  C \left(2^{\log (|\log t| +1) + 2}\right)^{3k_0}\nonumber\\
 &\leq&  C(4(1+|\log t|))^{3(\log (|\log t|+1)+2)}.
\end{eqnarray}

For the second factor of (\ref{eqfm}), we use Lemma \ref{sek} to
bound the denominator, together with Lemma \ref{sek0t} to bound
the numerator.  We get
\begin{eqnarray*}
 \prod_{k=k_0+1}^n \frac{\prob(SE_k(0,t) \ | \
        SE_{k-1}(0,t))}{\prob(SE_k(0) \ | \ SE_{k-1}(0))^2}
  &\leq &
    \prod_{k=k_0+1}^n
    \frac{(\beta_k \pi \sum_{2^k \leq j < 2^{k+1}}\frac{1}{j^3})^2 +\frac{C k}{2^{5k}}}
        {(\beta_k \pi \sum_{2^k \leq j < 2^{k+1}}\frac{1}{j^3} -
            \frac{C k}{2^{3k}})^2}\\
  & \leq &
    \prod_{k=k_0+1}^n \frac{(\beta_k \pi \sum_{2^k \leq j < 2^{k+1}}\frac{1}{j^3})^2 +\frac{C k}{2^{5k}}} {(\beta_k \pi
        \sum_{2^k \leq j < 2^{k+1}}\frac{1}{j^3})^2 -\frac{C k}{2^{5k}}}\\
  & \leq & \prod_{k=k_0+1}^n \left(1 + \frac{C k }{2^k}\right) \\
  & \leq & C.
\end{eqnarray*}
Multiplying the two estimates we deduce that
\begin{equation}\label{ft} f(t,M) \leq C(4(1+|\log t|))^{3(\log (|\log
t|+1)+2)}\end{equation} for any $M$. Since
$$\int_0^1 C(4(1+|\log t|))^{3(\log (|\log t|+1))+2} \, dt < \infty, $$
the proof is complete.
\end{pfof}

\section{Proof of Theorem \ref{expdim1}}
The proof of Theorem \ref{expdim1} follows from Theorem \ref{fmt}
and a second moment argument exactly as in \cite{H1}.  We give it
here for the sake of completeness.

\begin{pfof}{Theorem \ref{expdim1}}
Define $$ E_M(t) = \bigcap_{1\leq i\leq M}SE_i(t), $$

$$T_M=\{t:t \in [0,1] \mbox{ and }\esc_M(t) \mbox{ occurs}\}$$
and
$$T=\cap_1^{\infty}\overline{T_M}.$$
Now we show that $T$ is contained in the union of $\Sexc$ and the
countable set $$\Lambda=(\cup_{n,m}\tau_n^{(m)})\cup 1.$$ If $t
\in \cap_1^{\infty} T_M$ then $t \in \Sexc$. So if $t \in T
\setminus \Sexc$ then $t$ is contained in the boundary of $T_M$
for some $M$. For any $M$ the boundary of $T_M$ is contained in
$\Lambda$. Thus if $t \in T \setminus \Sexc$ then $t\in \Lambda$
and
$$T \subset \Sexc \cup \Lambda.$$
As $\Lambda$ is countable if $T$ has dimension one with positive
probability then so does $\Sexc$.

By Theorem \ref{fmt} there exists $f(t)$ such that
 $$\int_0^1 f(t) dt < \infty$$ and for all $M$
 \begin{equation} \label{downlow}
 \frac{\prob(\esc_M(0,t))}{(\prob(\esc_M(0)))^2}<f(M,t)< f(t).
 \end{equation}
Let $\leb(*)$ denote Lebesgue measure on $[0,1]$.  Then we get
\begin{eqnarray}
\expe(\leb(T_M)^2)
 & =    & \int_0^1 \int_0^1 \prob(\esc_M(r,s)) dr \times ds \label{fubini}\\
 & = & \int_0^1 \int_0^1 \prob(\esc_M(0,|s-r|)) dr \times ds \label{inv}\\
 &\leq  & \int_0^1 2\int_0^1 \prob(\esc_M(0,t)) dt \times ds \nonumber\\
 &\leq  & 2 \int_0^1 f(t)  \prob(\esc_M(0))^2 dt \label{downlower}\\
 &\leq  & 2  \prob(\esc_M(0))^2\int_0^1 f(t)  dt. \label{downlowest}
\end{eqnarray}
The equality  (\ref{fubini}) is true by Fubini's theorem,
(\ref{inv}) is true because
$$\esc_M(a,b)=\esc_M(b,a)=\esc_M(0,|b-a|)$$
and (\ref{downlower}) follows from (\ref{downlow}).


By Jensen's inequality if $h(x)=0$ for all $x \notin A$ then
\be \label{jensen}
\expe(h^2) \geq \frac{\expe (h)^2}{ \prob(A)}.
\ee
Then we get
\begin{eqnarray}
2  \prob(\esc_M(0))^2\int_0^1 f(t)  dt & \geq &
\expe(\leb(T_M)^2) \label{ref2}\\
 & \geq & \frac{\expe (\leb(T_M))^2}{\prob(T_M \neq \emptyset)}
     \label{ref6}\\
 & \geq & \frac{\prob(\esc_M(0))^2}{\prob(T_M \neq \emptyset)}, \nonumber
\end{eqnarray}
where (\ref{ref2}) is a restatement of (\ref{downlowest}), and
(\ref{ref6}) follows from (\ref{jensen}) with $\leb(T_M)$ and $T_M
\neq \emptyset$ in place of $h$ and $A$.

Thus for all $M$
$$\prob(T_M \neq \emptyset)
    \geq \frac{1}{2\int_{0}^{1}f(t) \ dt}>0. $$
As $T$ is the intersection of the nested sequence of compact sets
${\overline T_M}$
$$
\prob(T \neq \emptyset)
    = \lim_{M \to \infty} \prob({\overline T_M} \neq \emptyset)
    = \lim_{M \to \infty} \prob(T_M  \neq \emptyset) \geq
    \frac{1}{2\int_{0}^{1}f(t) \ dt }.$$

Now we show that the dimensions of $T$ and $\Sexc$ are one. By
Lemma 5.1 of \cite{peres} for any $\beta<1$ there exists a random
nested sequence of compact sets $F_{k} \subset [0,1]$ such that
\begin{equation} \label{star2}
\prob(r \in F_{k})\geq C(s_k)^{-\beta}
\end{equation}
and
\begin{equation}\label{star1}
\prob(r,t \in F_{k})\leq C(s_k)^{-2\beta}|r-t|^{-\beta}.
\end{equation}
These sets also have the property that for any set $T$ if
\begin{equation} \label{yuval}
\prob(T \cap (\cap_1^{\infty} F_k) \neq \emptyset)>0
\end{equation}
then $T$ has dimension at least $\beta$. We construct $F_k$ to be
independent of the dynamical random walk. So by (\ref{ft}),
(\ref{star2}) and  (\ref{star1}) we get
\begin{equation} \label{dim}
\frac{\prob(r,t \in T_M \cap F_M)}{\prob(r \in T_M \cap F_M)^2}
    \leq C(1+|\,\log |r-t|\,|)^{3|\log (|\,\log |r-t|\, |+1)|}|r-t|^{-\beta}.
\end{equation}

Since $$\int_0^1 C(1+|\log t|)^{3|\log (|\log t|+1)|}t^{-\beta} \
dt < \infty, $$ the same second moment argument as above and
(\ref{dim}) implies that with positive probability $T$ satisfies
(\ref{yuval}). Thus $T$ has dimension $\beta$ with positive
probability. As
$$T \subset (\Sexc \cap [0,1]) \cup \Lambda,$$
and $\Lambda$ is countable, the dimension of the set of $\Sexc
\cap [0,1]$ is at least $\beta$ with positive probability. By the
ergodic theorem the dimension of the set of $\Sexc$ is at least
$\beta$ with probability one. As this holds for all $\beta<1$ the
dimension of $\Sexc$ is one a.s.

\end{pfof}

\

\subsection{Extending the technique}
We end the paper with three remarks on the scope of the techniques
in this paper.

\begin{enumerate}
\item The techniques in this paper can be extended to show other
sets of exceptional times exist. For any two finite sets $E$ and
$V$, for which it is possible to visit every point in $V$ without
hitting $E$, there is a.s. an exceptional time $t$ such that each
point in $V$ is visited infinitely often, while the set $E$
is not visited at all. \\
\item Not every set can be avoided by the random walk, as shown in
the following claim:
\begin{claim}
Let $L$ be a subset of $\Z^2$, with the property that there exits
some $M>0$, such that every point in $\Z^2$ is at distance at most
$M$ from some point in $L$.  Then
\[ \prob(\exists t\in[0,1] \ : |\{n: \ S_n(t)\notin L\}| < \infty)
= 0\] i.e there are a.s. no exceptional times at which $L$ is
visited only finitely often.
\end{claim}
\begin{proof}
First, we can reduce to the case that there is some (possible) path
from $0$ that misses $L$. Otherwise we can just drop points out of
$L$ until such a path exists without ruining the desired condition
on $L$

Second, we note that it is enough to prove that there are a.s.\ no
times at which the random walk \textbf{never} hits $L$. This
follows from the Markovian property of random walk, and the fact
that changing finitely many moves is enough to make a walk that
 visits $L$ a finite number of times miss $L$ all together.

Let
$$Q_n = \{t\in[0,1]\ : \ S_i(t)\notin L \ \ \ 0\leq \forall i \leq n\}.$$
Then an exceptional time exists if and only if
$$\bigcap_{n\geq 0} Q_n \neq \emptyset.$$

 The condition on $L$ ensures that there is some $\epsilon > 0$, such
that for any simple random walk on $Z^2$, regardless of current
position or history, there is a probability of at least $\epsilon$
of hitting $L$ in the next $2M$ steps. Therefore there exists some
constants $A>0$ and $0 <c<1$ such that
$$\prob(0 \in Q_n) \leq Ac^n$$
 For any $i\geq 0$ define
 $$m_i = \min \{m\geq 0 : \tau_i^{(m)}>1\}.$$
Then
$$RE_i = \left\{\tau_i^{(j)}\right\}_{j=0}^{m_i}$$
is the set of times in $[0,1]$ in which the $i-th$ move gets
resampled. Let $H_n = \bigcup_{i=0}^n RE_i$. Then $H_n$ is the set
of all times at which one of the first $n$ moves of the walk are
re-sampled, and if we order $H_n$, then between two consecutive
times in $H_n$, the first $n$ steps of the walk do not change.
Thus $Q_n \neq \emptyset$ if and only if there exists some $\tau
\in H_n$ for which $\tau \in Q_n$. Since $\{S_n(t)\}$ behaves like
a simple random walk for all $t$, we have, by linearity of
expectation
\[\Ex(\# \tau\in H_n \ : \ \tau\in Q_n) = \Ex(|H_n|)\prob(0\in
Q_n) \leq 2Anc^{-n}\] Thus
\[\prob(Q_n \neq \emptyset) \leq \Ex(\# \tau\in H_n \ : \
\tau\in Q_n) \leq 2Anc^{-n}\]

And consequently $$\prob\left(\bigcap_{n\geq 0} Q_n \neq
\emptyset\right) = 0.
$$
So no exceptional times exist.
\end{proof}

 \item We finish with an open question: \\
 Is there a set $A \subset \Z^2$ such that both $A$ and $A^c$ are
 infinite and almost surely there
 are exceptional times in which every element of $A$ is hit finitely often and every element of
 $A^c$ is hit infinitely often?
 \end{enumerate}

\end{document}